

\baselineskip=14pt
\parskip=10pt

\magnification=\magstephalf

\def\1{{\overline{1}}}
\def\2{{\overline{2}}}
\parindent=0pt
\overfullrule=0in

\def\frac#1#2{{#1 \over #2}}
\centerline
{\bf
Searching for Disjoint Covering Systems with Precisely One Repeated Modulus
}

\bigskip
\centerline
{\it By Shalosh B. EKHAD, Aviezri S. FRAENKEL, and Doron ZEILBERGER}

{\bf Introduction}

A set of arithmetical sequences
$$
a_1\, (\bmod{ \,\, m_1}) \quad, \quad  a_2 \, (\bmod{\,\, m_2}) \quad, \quad \dots \quad , \quad  a_k \, (\bmod{\,\,m_k}) \quad \quad ,
$$
with 
$$
m_1 \leq m_2 \leq \dots \leq m_k \quad \quad ,
$$
is called a {\it disjoint covering system} (alias  {\it exact covering system})
if every positive integer belongs to {\bf exactly} one of the sequences.

Mirski, Newman, Davenport and Rado (See [3])
famously proved  that the moduli can't all be distinct. In fact the two largest moduli must be
equal, i.e.
$$
m_{k-1}=m_k \quad .
$$

This raises the natural question: 

{\it How close can you get to getting distinct moduli?}

Of course, for any positive integer $r$, the system
$$
0 \, (\bmod{\,\,2}) \quad, \quad
1 \, (\bmod{\,\,4}) \quad, \quad
3 \, (\bmod{\,\, 8}) \quad, \quad 
\dots \quad , \quad
$$
$$
2^{r-2} -1 \, (\bmod{2^{r-1}}) \quad, \quad 
2^{r-1} -1 \, (\bmod{2^r}) \quad, \quad 
2^{r} -1 \, (\bmod{2^r}) \quad, 
$$
has the property that all the moduli are distinct except for the largest modulus, $2^r$,  that only appears twice.
So there are {\it infinitely} many systems where the moduli are all distinct except the largest moduli 
repeating twice. We label these as {\it trivial}.

More generally, for any given positive integer $r$, 
once  you have one system  with all distinct moduli, except for the largest, that is repeated, $r$ times,
one can manufacture infinitely many examples, by the "add-2" process. Multiply all the
moduli by $2$ and prefix $2$ at the beginning.

But if you insist that the smallest modulus is at least three, then, for any given $r$, there are 
(conjecturally, but most probably) only finitely many such systems. This raises the
natural question of finding all of them for small $r$.

This was done by Berger, Felzenbaum, and Fraenkel([2]), for $3 \leq r \leq 9$,  and extended by
Zeleke and Simpson ([3]) for $10 \leq r \leq 12$ ([4]). In this article we continue the list all
the way to $r=32$.

We must confess that, while we have proved that all the systems that we have found are indeed correct, we did not
bother to prove that all our lists are complete. What we do know, {\it for sure},  is that we found all such systems where
the largest modulus is $\leq 600$.
Nevertheless we strongly believe, on heuristic and experimental grounds, that there aren't any new ones.

The results below were found using the Maple package {\tt BFF.txt}, available, free of charge, from the
following url:

{\tt http://www.math.rutgers.edu/\~{}zeilberg/tokhniot/BFF.txt} \quad .

The front of the present article

{\tt http://www.math.rutgers.edu/\~{}zeilberg/mamarim/mamarimhtml/ecs.html} \quad ,
 
contains links to sample input and output files, that the readers are welcome to extend, if they have
computer time to spare.

{\bf The output: the list of all  Exact Covering Systems whose smallest modulus is at least 3 and
whose moduli are all distinct except the largest one that is repeated at most 32 times}

We exclude from the list the trivial systems consisting of only one modulus.

                                                                      {\bf 2  Repeats}

                                                                                None

                                                                       {\bf 3  Repeats}

                                                                                None

                                                                       {\bf 4  Repeats}

                                                               There  is 1 system here it is:

                                                                        3 6 ( 4  times)

                                                                            {\bf 5  Repeats}

                                                                                None

                                                                            {\bf 6  Repeats}

                                                               There are 3 systems here they are:

                                                                        3 9 ( 6  times)

                                                                        4 8 ( 6  times)

                                                                      3 6 12 ( 6  times)

                                                                            {\bf 7  Repeats}

                                                               There are 2 systems here they are:

                                                                      4 6 12 ( 7  times)

                                                                     3 6 9 18 ( 7  times)

                                                                            {\bf 8  Repeats}

                                                               There are 2 systems here they are:

                                                                        3 12 ( 8  times)

                                                                        5 10 ( 8  times)

                                                                            {\bf 9  Repeats}

                                                               There are 4 systems here they are:

                                                                        4 12 ( 9  times)

                                                                      3 6 18 ( 9  times)

                                                                   4 6 8 12 24 ( 9  times)

                                                                 3 6 9 12 18 36 ( 9  times)

                                                                            {\bf 10  Repeats}

                                                               There are 5 systems here they are:

                                                                       3 15 ( 10  times)

                                                                       6 12 ( 10  times)

                                                                      3 9 18 ( 10  times)

                                                                      4 8 16 ( 10  times)

                                                                    3 6 12 24 ( 10  times)

                                                                            {\bf 11  Repeats}

                                                               There are 2 systems here they are:

                                                                    4 6 8 24 ( 11  times)

                                                                  3 6 9 12 36 ( 11  times)

                                                                            {\bf 12  Repeats}

                                                               There are 7 systems here they are:

                                                                       3 18 ( 12  times)

                                                                       4 16 ( 12  times)

                                                                       5 15 ( 12  times)

                                                                       7 14 ( 12  times)

                                                                      3 6 24 ( 12  times)

                                                                    4 6 12 24 ( 12  times)

                                                                  3 6 9 18 36 ( 12  times)

                                                                            {\bf 13  Repeats}

                                                               There are 7 systems here they are:

                                                                     4 10 20 ( 13  times)

                                                                      6 9 18 ( 13  times)

                                                                    3 6 15 30 ( 13  times)

                                                                    4 8 12 24 ( 13  times)

                                                                  3 6 12 18 36 ( 13  times)

                                                              4 6 8 12 16 24 48 ( 13  times)

                                                            3 6 9 12 18 24 36 72 ( 13  times)

                                                                            {\bf 14  Repeats}

                                                               There are 6 systems here they are:

                                                                       3 21 ( 14  times)

                                                                       8 16 ( 14  times)

                                                                     3 12 24 ( 14  times)

                                                                      4 6 24 ( 14  times)

                                                                     5 10 20 ( 14  times)

                                                                    3 6 9 36 ( 14  times)

                                                                            {\bf 15  Repeats}

                                                               There are 10 systems here they are:

                                                                       4 20 ( 15  times)

                                                                       6 18 ( 15  times)

                                                                      3 6 30 ( 15  times)

                                                                      3 9 27 ( 15  times)

                                                                      4 8 24 ( 15  times)

                                                                    3 6 12 36 ( 15  times)

                                                                    6 8 12 24 ( 15  times)

                                                                  3 9 12 18 36 ( 15  times)

                                                                4 6 8 12 16 48 ( 15  times)

                                                              3 6 9 12 18 24 72 ( 15  times)

                                                                            {\bf 16  Repeats}

                                                               There are 9 systems here they are:

                                                                       3 24 ( 16  times)

                                                                       5 20 ( 16  times)

                                                                       9 18 ( 16  times)

                                                                     4 12 24 ( 16  times)

                                                                    3 6 18 36 ( 16  times)

                                                                  4 6 12 18 36 ( 16  times)

                                                                3 6 9 18 27 54 ( 16  times)

                                                                4 6 8 12 24 48 ( 16  times)

                                                              3 6 9 12 18 36 72 ( 16  times)

                                                                            {\bf 17  Repeats}

                                                               There are 4 systems here they are:

                                                                      6 8 24 ( 17  times)

                                                                    3 9 12 36 ( 17  times)

                                                                4 6 8 16 24 48 ( 17  times)

                                                              3 6 9 12 24 36 72 ( 17  times)

                                                                            {\bf 18  Repeats}

                                                               There are 16 systems here they are:

                                                                       3 27 ( 18  times)

                                                                       4 24 ( 18  times)

                                                                       7 21 ( 18  times)

                                                                       10 20 ( 18  times)

                                                                      3 6 36 ( 18  times)

                                                                     3 15 30 ( 18  times)

                                                                     6 12 24 ( 18  times)

                                                                    3 9 18 36 ( 18  times)

                                                                    4 6 12 36 ( 18  times)

                                                                    4 8 16 32 ( 18  times)

                                                                  3 6 9 18 54 ( 18  times)

                                                                  3 6 12 24 48 ( 18  times)

                                                                  4 6 8 12 48 ( 18  times)

                                                                3 6 9 12 18 72 ( 18  times)

                                                            4 6 8 12 18 24 36 72 ( 18  times)

                                                          3 6 9 12 18 27 36 54 108 ( 18  times)

                                                                            {\bf 19  Repeats}

                                                               There are 13 systems here they are:

                                                                     4 14 28 ( 19  times)

                                                                     8 12 24 ( 19  times)

                                                                    3 6 21 42 ( 19  times)

                                                                   3 12 18 36 ( 19  times)

                                                                    4 6 18 36 ( 19  times)

                                                                   5 10 15 30 ( 19  times)

                                                                  3 6 9 27 54 ( 19  times)

                                                                  4 6 8 16 48 ( 19  times)

                                                                  4 8 10 20 40 ( 19  times)

                                                                3 6 9 12 24 72 ( 19  times)

                                                                3 6 12 15 30 60 ( 19  times)

                                                                4 6 12 16 24 48 ( 19  times)

                                                              3 6 9 18 24 36 72 ( 19  times)

                                                                            {\bf 20  Repeats}

                                                               There are 12 systems here they are:

                                                                       3 30 ( 20  times)

                                                                       5 25 ( 20  times)

                                                                       6 24 ( 20  times)

                                                                       11 22 ( 20  times)

                                                                      3 9 36 ( 20  times)

                                                                      4 8 32 ( 20  times)

                                                                    3 6 12 48 ( 20  times)

                                                                   6 10 15 30 ( 20  times)

                                                                  4 6 8 24 48 ( 20  times)

                                                                3 6 9 12 36 72 ( 20  times)

                                                              4 6 8 12 18 24 72 ( 20  times)

                                                            3 6 9 12 18 27 36 108 ( 20  times)

                                                                            {\bf 21  Repeats}

                                                               There are 18 systems here they are:

                                                                       4 28 ( 21  times)

                                                                       8 24 ( 21  times)

                                                                      3 6 42 ( 21  times)

                                                                     3 12 36 ( 21  times)

                                                                      4 6 36 ( 21  times)

                                                                     5 10 30 ( 21  times)

                                                                    3 6 9 54 ( 21  times)

                                                                    4 8 10 40 ( 21  times)

                                                                  3 6 12 15 60 ( 21  times)

                                                                  4 6 12 16 48 ( 21  times)

                                                                  6 9 12 18 36 ( 21  times)

                                                                3 6 9 18 24 72 ( 21  times)

                                                                4 8 12 16 24 48 ( 21  times)

                                                              3 6 12 18 24 36 72 ( 21  times)

                                                              4 6 8 12 18 36 72 ( 21  times)

                                                            3 6 9 12 18 27 54 108 ( 21  times)

                                                          4 6 8 12 16 24 32 48 96 ( 21  times)

                                                        3 6 9 12 18 24 36 48 72 144 ( 21  times)

\vfill\eject

                                                                            {\bf 22  Repeats}

                                                               There are 19 systems here they are:

                                                                       3 33 ( 22  times)

                                                                       12 24 ( 22  times)

                                                                     3 18 36 ( 22  times)

                                                                     4 16 32 ( 22  times)

                                                                     5 15 30 ( 22  times)

                                                                     6 10 30 ( 22  times)

                                                                     7 14 28 ( 22  times)

                                                                    3 6 24 48 ( 22  times)

                                                                    3 9 15 45 ( 22  times)

                                                                    4 6 8 48 ( 22  times)

                                                                   4 12 18 36 ( 22  times)

                                                                  3 6 9 12 72 ( 22  times)

                                                                  3 6 18 27 54 ( 22  times)

                                                                  4 6 12 24 48 ( 22  times)

                                                                3 6 9 18 36 72 ( 22  times)

                                                              4 6 8 12 24 36 72 ( 22  times)

                                                            3 6 9 12 18 36 54 108 ( 22  times)

                                                      4 6 8 12 16 18 24 36 48 72 144 ( 22  times)

                                                   3 6 9 12 18 24 27 36 54 72 108 216 ( 22  times)

                                                                            {\bf 23  Repeats}

                                                               There are 12 systems here they are:

                                                                     6 15 30 ( 23  times)

                                                                    4 8 20 40 ( 23  times)

                                                                    6 9 12 36 ( 23  times)

                                                                  3 6 12 30 60 ( 23  times)

                                                                  4 6 16 24 48 ( 23  times)

                                                                  4 8 12 16 48 ( 23  times)

                                                                3 6 9 24 36 72 ( 23  times)

                                                                3 6 12 18 24 72 ( 23  times)

                                                                4 6 8 12 18 72 ( 23  times)

                                                              3 6 9 12 18 27 108 ( 23  times)

                                                            4 6 8 12 16 24 32 96 ( 23  times)

                                                          3 6 9 12 18 24 36 48 144 ( 23  times)

                                                                            {\bf 24  Repeats}

                                                               There are 24 systems here they are:

                                                                       3 36 ( 24  times)

                                                                       4 32 ( 24  times)

                                                                       5 30 ( 24  times)

                                                                       7 28 ( 24  times)

                                                                       9 27 ( 24  times)

                                                                       13 26 ( 24  times)

                                                                      3 6 48 ( 24  times)

                                                                     4 12 36 ( 24  times)

                                                                    3 6 18 54 ( 24  times)

                                                                    4 6 12 48 ( 24  times)

                                                                   4 10 20 40 ( 24  times)

                                                                    6 9 18 36 ( 24  times)

                                                                  3 6 9 18 72 ( 24  times)

                                                                  3 6 15 30 60 ( 24  times)

                                                                  4 8 12 24 48 ( 24  times)

                                                                3 6 12 18 36 72 ( 24  times)

                                                                4 6 8 12 24 72 ( 24  times)

                                                              3 6 9 12 18 36 108 ( 24  times)

                                                              4 6 8 18 24 36 72 ( 24  times)

                                                            3 6 9 12 27 36 54 108 ( 24  times)

                                                            4 6 8 12 16 24 48 96 ( 24  times)

                                                          3 6 9 12 18 24 36 72 144 ( 24  times)

                                                        4 6 8 12 16 18 24 36 48 144 ( 24  times)

                                                      3 6 9 12 18 24 27 36 54 72 216 ( 24  times)

                                                                            {\bf 25  Repeats}

                                                               There are 23 systems here they are:

                                                                       6 30 ( 25  times)

                                                                      3 9 45 ( 25  times)

                                                                      4 8 40 ( 25  times)

                                                                     4 18 36 ( 25  times)

                                                                     10 15 30 ( 25  times)

                                                                    3 6 12 60 ( 25  times)

                                                                    3 6 27 54 ( 25  times)

                                                                    4 6 16 48 ( 25  times)

                                                                   6 12 18 36 ( 25  times)

                                                                  3 6 9 24 72 ( 25  times)

                                                                  3 9 18 27 54 ( 25  times)

                                                                  4 8 16 24 48 ( 25  times)

                                                                3 6 12 24 36 72 ( 25  times)

                                                                4 6 8 12 36 72 ( 25  times)

                                                                4 6 12 20 30 60 ( 25  times)

                                                                6 8 12 16 24 48 ( 25  times)

                                                              3 6 9 12 18 54 108 ( 25  times)

                                                              3 6 9 18 30 45 90 ( 25  times)

                                                              3 9 12 18 24 36 72 ( 25  times)

                                                            4 6 8 12 16 32 48 96 ( 25  times)

                                                          3 6 9 12 18 24 48 72 144 ( 25  times)

                                                        4 6 8 12 16 18 24 36 72 144 ( 25  times)

                                                     3 6 9 12 18 24 27 36 54 108 216 ( 25  times)

                                                                            {\bf 26  Repeats}

                                                               There are 19 systems here they are:

                                                                       3 39 ( 26  times)

                                                                       14 28 ( 26  times)

                                                                     3 21 42 ( 26  times)

                                                                     4 10 40 ( 26  times)

                                                                      6 9 36 ( 26  times)

                                                                     8 16 32 ( 26  times)

                                                                    3 6 15 60 ( 26  times)

                                                                   3 12 24 48 ( 26  times)

                                                                    4 6 24 48 ( 26  times)

                                                                    4 8 12 48 ( 26  times)

                                                                   5 10 20 40 ( 26  times)

                                                                  3 6 9 36 72 ( 26  times)

                                                                  3 6 12 18 72 ( 26  times)

                                                                4 6 8 18 24 72 ( 26  times)

                                                              3 6 9 12 27 36 108 ( 26  times)

                                                              4 6 8 12 16 24 96 ( 26  times)

                                                            3 6 9 12 18 24 36 144 ( 26  times)

                                                        4 6 8 12 16 18 24 48 72 144 ( 26  times)

                                                     3 6 9 12 18 24 27 36 72 108 216 ( 26  times)

                                                                            {\bf 27  Repeats}

                                                               There are 27 systems here they are:

                                                                       4 36 ( 27  times)

                                                                       10 30 ( 27  times)

                                                                      3 6 54 ( 27  times)

                                                                     3 15 45 ( 27  times)

                                                                     6 12 36 ( 27  times)

                                                                    3 9 18 54 ( 27  times)

                                                                    4 8 16 48 ( 27  times)

                                                                   9 12 18 36 ( 27  times)

                                                                  3 6 12 24 72 ( 27  times)

                                                                  4 6 8 12 72 ( 27  times)

                                                                  4 6 12 20 60 ( 27  times)

                                                                 4 12 16 24 48 ( 27  times)

                                                                  6 8 12 16 48 ( 27  times)

                                                                3 6 9 12 18 108 ( 27  times)

                                                                3 6 9 18 30 90 ( 27  times)

                                                                3 6 18 24 36 72 ( 27  times)

                                                                3 9 12 18 24 72 ( 27  times)

                                                                4 6 8 18 36 72 ( 27  times)

                                                              3 6 9 12 27 54 108 ( 27  times)

                                                              4 6 8 12 16 32 96 ( 27  times)

                                                              4 6 12 18 24 36 72 ( 27  times)

                                                            3 6 9 12 18 24 48 144 ( 27  times)

                                                            3 6 9 18 27 36 54 108 ( 27  times)

                                                            4 6 8 12 24 32 48 96 ( 27  times)

                                                          3 6 9 12 18 36 48 72 144 ( 27  times)

                                                          4 6 8 12 16 18 24 36 144 ( 27  times)

                                                        3 6 9 12 18 24 27 36 54 216 ( 27  times)

                                                                            {\bf 28  Repeats}

                                                               There are 26 systems here they are:

                                                                       3 42 ( 28  times)

                                                                       5 35 ( 28  times)

                                                                       8 32 ( 28  times)

                                                                       15 30 ( 28  times)

                                                                     3 12 48 ( 28  times)

                                                                      4 6 48 ( 28  times)

                                                                     4 20 40 ( 28  times)

                                                                     5 10 40 ( 28  times)

                                                                     6 18 36 ( 28  times)

                                                                    3 6 9 72 ( 28  times)

                                                                    3 6 30 60 ( 28  times)

                                                                    3 9 27 54 ( 28  times)

                                                                    4 8 24 48 ( 28  times)

                                                                  3 6 12 36 72 ( 28  times)

                                                                  4 6 12 30 60 ( 28  times)

                                                                  6 8 12 24 48 ( 28  times)

                                                                3 6 9 18 45 90 ( 28  times)

                                                                3 9 12 18 36 72 ( 28  times)

                                                                4 6 8 24 36 72 ( 28  times)

                                                              3 6 9 12 36 54 108 ( 28  times)

                                                              4 6 8 12 16 48 96 ( 28  times)

                                                            3 6 9 12 18 24 72 144 ( 28  times)

                                                          4 6 8 12 16 18 24 48 144 ( 28  times)

                                                        3 6 9 12 18 24 27 36 72 216 ( 28  times)

                                                        4 6 8 12 16 18 36 48 72 144 ( 28  times)

                                                     3 6 9 12 18 24 27 54 72 108 216 ( 28  times)

                                                                            {\bf 29  Repeats}

                                                               There are 21 systems here they are:

                                                                     9 12 36 ( 29  times)

                                                                   4 12 16 48 ( 29  times)

                                                                   8 10 20 40 ( 29  times)

                                                                  3 6 18 24 72 ( 29  times)

                                                                 3 12 15 30 60 ( 29  times)

                                                                  4 6 8 18 72 ( 29  times)

                                                                  4 8 14 28 56 ( 29  times)

                                                                  6 8 16 24 48 ( 29  times)

                                                                3 6 9 12 27 108 ( 29  times)

                                                                3 6 12 21 42 84 ( 29  times)

                                                                3 9 12 24 36 72 ( 29  times)

                                                                4 6 12 18 24 72 ( 29  times)

                                                               4 10 12 20 30 60 ( 29  times)

                                                              3 6 9 18 27 36 108 ( 29  times)

                                                              3 6 15 18 30 45 90 ( 29  times)

                                                              4 6 8 12 24 32 96 ( 29  times)

                                                            3 6 9 12 18 36 48 144 ( 29  times)

                                                            4 6 8 16 24 32 48 96 ( 29  times)

                                                          3 6 9 12 24 36 48 72 144 ( 29  times)

                                                          4 6 8 12 16 18 24 72 144 ( 29  times)

                                                       3 6 9 12 18 24 27 36 108 216 ( 29  times)

                                                                            {\bf 30  Repeats}

                                                               There are 39 systems here they are:

                                                                       3 45 ( 30  times)

                                                                       4 40 ( 30  times)

                                                                       6 36 ( 30  times)

                                                                       7 35 ( 30  times)

                                                                       11 33 ( 30  times)

                                                                       16 32 ( 30  times)

                                                                      3 6 60 ( 30  times)

                                                                      3 9 54 ( 30  times)

                                                                     3 24 48 ( 30  times)

                                                                      4 8 48 ( 30  times)

                                                                     5 20 40 ( 30  times)

                                                                     9 18 36 ( 30  times)

                                                                    3 6 12 72 ( 30  times)

                                                                    4 6 12 60 ( 30  times)

                                                                   4 12 24 48 ( 30  times)

                                                                    6 8 12 48 ( 30  times)

                                                                   6 14 21 42 ( 30  times)

                                                                  3 6 9 18 90 ( 30  times)

                                                                  3 6 18 36 72 ( 30  times)

                                                                  3 9 12 18 72 ( 30  times)

                                                                  4 6 8 24 72 ( 30  times)

                                                                  4 6 20 30 60 ( 30  times)

                                                                3 6 9 12 36 108 ( 30  times)

                                                                3 6 9 30 45 90 ( 30  times)

                                                                4 6 8 12 16 96 ( 30  times)

                                                                4 6 12 18 36 72 ( 30  times)

                                                              3 6 9 12 18 24 144 ( 30  times)

                                                              3 6 9 18 27 54 108 ( 30  times)

                                                              4 6 8 12 24 48 96 ( 30  times)

                                                              4 8 12 18 24 36 72 ( 30  times)

                                                             6 10 12 15 20 30 60 ( 30  times)

                                                            3 6 9 12 18 36 72 144 ( 30  times)

                                                           3 6 12 18 27 36 54 108 ( 30  times)

                                                          4 6 8 12 16 18 36 48 144 ( 30  times)

                                                        3 6 9 12 18 24 27 54 72 216 ( 30  times)

                                                        4 6 8 12 16 24 36 48 72 144 ( 30  times)

                                                     3 6 9 12 18 24 36 54 72 108 216 ( 30  times)

                                               4 6 8 12 16 18 24 32 36 48 72 96 144 288 ( 30  times)

                                            3 6 9 12 18 24 27 36 48 54 72 108 144 216 432 ( 30  times)

                                                                            {\bf 31  Repeats}

                                                               There are 35 systems here they are:

                                                                     4 22 44 ( 31  times)

                                                                     8 10 40 ( 31  times)

                                                                     12 18 36 ( 31  times)

                                                                    3 6 33 66 ( 31  times)

                                                                   3 12 15 60 ( 31  times)

                                                                   3 18 27 54 ( 31  times)

                                                                    4 8 14 56 ( 31  times)

                                                                   4 16 24 48 ( 31  times)

                                                                    6 8 16 48 ( 31  times)

                                                                   7 14 21 42 ( 31  times)

                                                                  3 6 12 21 84 ( 31  times)

                                                                  3 6 24 36 72 ( 31  times)

                                                                  3 9 12 24 72 ( 31  times)

                                                                  4 6 8 36 72 ( 31  times)

                                                                 4 10 12 20 60 ( 31  times)

                                                                 6 12 16 24 48 ( 31  times)

                                                                3 6 9 12 54 108 ( 31  times)

                                                                3 6 15 18 30 90 ( 31  times)

                                                                3 9 18 24 36 72 ( 31  times)

                                                                4 6 12 24 36 72 ( 31  times)

                                                              3 6 9 18 36 54 108 ( 31  times)

                                                              4 6 8 12 32 48 96 ( 31  times)

                                                              4 6 8 16 24 32 96 ( 31  times)

                                                              4 8 10 16 20 40 80 ( 31  times)

                                                            3 6 9 12 18 48 72 144 ( 31  times)

                                                            3 6 9 12 24 36 48 144 ( 31  times)

                                                           3 6 12 15 24 30 60 120 ( 31  times)

                                                            4 6 8 12 16 18 24 144 ( 31  times)

                                                          3 6 9 12 18 24 27 36 216 ( 31  times)

                                                          4 6 8 12 16 18 36 72 144 ( 31  times)

                                                          4 6 8 12 24 30 40 60 120 ( 31  times)

                                                       3 6 9 12 18 24 27 54 108 216 ( 31  times)

                                                        3 6 9 12 18 36 45 60 90 180 ( 31  times)

                                                        4 6 8 12 18 24 36 48 72 144 ( 31  times)

                                                     3 6 9 12 18 27 36 54 72 108 216 ( 31  times)

                                                                            {\bf 32  Repeats}

                                                               There are 29 systems here they are:

                                                                       3 48 ( 32  times)

                                                                       5 40 ( 32  times)

                                                                       9 36 ( 32  times)

                                                                       17 34 ( 32  times)

                                                                     4 12 48 ( 32  times)

                                                                     6 14 42 ( 32  times)

                                                                    3 6 18 72 ( 32  times)

                                                                    3 9 21 63 ( 32  times)

                                                                    4 6 20 60 ( 32  times)

                                                                    6 8 24 48 ( 32  times)

                                                                  3 6 9 30 90 ( 32  times)

                                                                  3 9 12 36 72 ( 32  times)

                                                                  4 6 12 18 72 ( 32  times)

                                                                 4 10 12 30 60 ( 32  times)

                                                                3 6 9 18 27 108 ( 32  times)

                                                                3 6 15 18 45 90 ( 32  times)

                                                                4 6 8 12 24 96 ( 32  times)

                                                                4 8 12 18 24 72 ( 32  times)

                                                               6 10 12 15 20 60 ( 32  times)

                                                              3 6 9 12 18 36 144 ( 32  times)

                                                             3 6 12 18 27 36 108 ( 32  times)

                                                              4 6 8 16 24 48 96 ( 32  times)

                                                            3 6 9 12 24 36 72 144 ( 32  times)

                                                          4 6 8 12 16 18 48 72 144 ( 32  times)

                                                          4 6 8 12 16 24 36 48 144 ( 32  times)

                                                       3 6 9 12 18 24 27 72 108 216 ( 32  times)

                                                        3 6 9 12 18 24 36 54 72 216 ( 32  times)

                                                  4 6 8 12 16 18 24 32 36 48 72 96 288 ( 32  times)

                                               3 6 9 12 18 24 27 36 48 54 72 108 144 432 ( 32  times)

\vfill\eject

{\bf Final Comments}

{\bf 1.} For a small number of repeats, the present computations and those done previously, 
confirm the proved and conjectured structures of such systems with a single repeated modulus. 
It is an open problem to characterize the structure for an arbitrary number of repeats.

{\bf 2.} Another direction of investigation would be to compute the number of repeats when we have exactly two distinct 
repeated moduli. A lower bound on the number of repeats for any number of repeated moduli is given in [1].
This direction may be conducive for solving {\bf 1}.

{\bf 3.} This article could have been written in many fewer pages (and had we published it in the
obsolete medium of print-journal, it would have been a good idea, in order to save trees), by
presenting the data concisely in a table, where the
values of $r$ appear in the leftmost column. Thus, for example, the table
entry for $r=14$ would be:
 
14(6): 3, {\bf 21}; 8, {\bf 16}; 3,12, {\bf 24}; 4,6, {\bf 24}; 5,10, {\bf 20}; 3,6, 9, {\bf 36} \quad .
 
(The parenthesized number next to $r$ is the number of systems. The number in boldface is the modulus that's repeated $r$ times.)

But for the sake of {\it readability} and {\it clarity}, and with the hope that no one of our readers will be stupid
enough to print out this paper, we decided to keep it in the present form.

{\bf References}

1. M.A. Berger, A. Felzenbaum, and A. S. Fraenkel, 
{\it Improvements to the Newman-Znam result for disjoint covering systems}, 
Acta Arithmetica {\bf 50}(1988), 1-13.

2. M.A. Berger, A. Felzenbaum, and A. S. Fraenkel, {\it Disjoint covering systems with precisely one multiple modulus}, 
Acta Arithmetica {\bf 50}(1988), 171-182.

3. P. Erd\H{o}s, On a problem concerning congruence systems. (Hungarian, with English summary) Mat. Lapok {\bf 3}(1952), 122-128.

4. M. Zeleke and J. Simpson, {\it On Disjoint Covering Systems with precisely one repeated modulus},
Advances in Applied Mathematics {\bf 23}(1999), 322-332.

\vfill\eject

\bigskip
\hrule

Shalosh B. Ekhad, c/o Doron Zeilberger, Department of Mathematics, Rutgers University (New Brunswick), Hill Center-Busch Campus, 110 Frelinghuysen
Rd., Piscataway, NJ 08854-8019, USA.  \hfill\break
ShaloshBEkhad at gmail dot com \quad .

\smallskip
Aviezri  S. Fraenkel, Department of Computer Science and Applied  Mathematics, Weizmann Institute of Science, Rehovot 76100, Israel ;
aviezri dot fraenkel at weizmann dot ac dot il \quad .
\smallskip
Doron Zeilberger, Department of Mathematics, Rutgers University (New Brunswick), Hill Center-Busch Campus, 110 Frelinghuysen
Rd., Piscataway, NJ 08854-8019, USA.  \hfill\break
zeilberg at math dot rutgers dot edu \quad .

\smallskip
\hrule
\medskip

{\bf Exclusively published in the Personal Journal of Shalosh B. Ekhad and Doron Zeilberger and arxiv.org, Nov. 12, 2015.}
\end